\documentclass[smallextended]{article}
\usepackage[english]{babel}
\usepackage[dvips]{graphicx}
\usepackage{amsthm}
\usepackage{exscale}
\usepackage{amsfonts}
\usepackage{amssymb}
\usepackage[intlimits,sumlimits]{amsmath}
\usepackage{amsmath,cases}
\usepackage[latin1]{inputenc}
\usepackage{mathrsfs}
\usepackage{fancyhdr}
\usepackage[misc,geometry]{ifsym}

\theoremstyle{definition}
\newtheorem{thm}{Theorem}

\newtheorem{prop}[thm]{Proposition}

\newtheorem{remark}{Remark}

\usepackage{a4wide}

\newcommand{\ud}{\mathrm{d}}

\newcommand{\xib}{\tilde{\boldsymbol{\xi}}}

\newcommand{\xisun}{\tilde{\boldsymbol{\xi}}^{(n)}}
\newcommand{\xisunm}{\tilde{\boldsymbol{\xi}}_{(n,m)}}

\newcommand{\inteone}{\int_{\mathbb{R}}}
\newcommand{\ind}{1 \! \textrm{l}}

\def\rd{\mathbb{R}^d}

\def\naturals{\mathbb{N}}

\def\reals{\mathbb{R}}

\newcommand{\pp}{\textsf{P}}
\newcommand{\ee}{\textsf{E}}
\newcommand{\probabilityspace}{(\Omega, \mathscr{F}, \textsf{P})}
\newcommand{\samplespace}{\mathbb{X}}

\newcommand{\pms}{[\samplespace]}
\newcommand{\randommeasure}{\tilde{\mathfrak{p}}}
\newcommand{\empiric}{\tilde{\mathfrak{e}}_n}
\newcommand{\pfrak}{\mathfrak{p}}
\newcommand{\Space}{\mathbb{S}}

\author{Donato Michele Cifarelli, Emanuele Dolera,  Eugenio Regazzini}

\title{Frequentistic approximations to Bayesian prevision \footnote{The term \emph{prevision} is a translation of the Italian \emph{previsione}, adopted by de Finetti, which is used in the English translation of his treatise (1970). See Translator's note on page 21.}\ \ \  of exchangeable random elements}

\begin{document}
\maketitle

\small{\textbf{Abstract} Given a sequence $\tilde{\xi}_1, \tilde{\xi}_2, \dots$ of $\samplespace$-valued, exchangeable random elements, let $q(\xisun)$ and $p_m(\xisun)$ stand for posterior and predictive distribution, respectively, given $\xisun := (\tilde{\xi}_1, \dots, \tilde{\xi}_n)$. We provide an upper bound for $\limsup_{n \rightarrow +\infty} b_n \ud_{[\pms]}(q(\xisun), \delta_{\empiric})$ and $\limsup_{n \rightarrow +\infty} b_n \ud_{[\samplespace^m]}(p_m(\xisun), \empiric^m)$, where $\empiric := \frac{1}{n} \sum_{i=1}^n \delta_{\tilde{\xi}_i}$ is the empirical measure, $\{b_n\}_{n \geq 1}$ is a suitable sequence of positive numbers increasing to $+\infty$, $\ud_{[\pms]}$ and $\ud_{[\samplespace^m]}$ denote distinguished weak probability distances on $[\pms]$ and $[\samplespace^m]$, respectively, with the proviso that $[\mathbb{S}]$ denotes the space of all probability measures on $\mathbb{S}$. A characteristic feature of our work is that the aforesaid bounds are established under the law of the $\tilde{\xi}_n$'s, unlike the more common literature on Bayesian consistency, where they are studied with respect to product measures $\pfrak_0^{\infty}$, as $\pfrak_0$ varies among the admissible determinations of a random probability measure.}

\section{Introduction} \label{sect:intro}

In the present paper the term \emph{prevision} will be used to designate both any activity directed to evaluation of probabilities of future (or, at least, till not known) events on the basis of an observed frequency, and the result of such an activity. Thus, prevision mingles with probabilistic inductive reasoning, and an important field of application of prevision is that of statistical problems, classically characterized by the circumstance that the events considered therein are generally thought of as \emph{analogous} events. Frequentistic approaches to statistics look at observable single events---or more general random elements $\tilde{\xi}_1, \tilde{\xi}_2, \dots$ taking values in some space $\samplespace$, like in the rest of the present work---as independent and identically distributed (i.i.d.) with a common probability distribution (p.d.) that can be approximated by observed (empirical) frequency. Laws of large numbers and allied results are then invoked to assert that such an approximation improves as the number of observations goes to infinity. Bayesian statisticians translate the aforesaid analogy into a less restrictive property, that is the \emph{exchangeability} of the $\tilde{\xi}_n$'s. As a consequence, any correct expression of Bayesian prevision must rely on a conditional p.d. for till now unknown observable random elements, given the frequency distribution of observed random elements. The expectation, due to the analogy of the observable elements here realized in the form of exchangeability, is that we are willing to be influenced more and more by the observed frequency as the size of experience increases. The present paper hinges upon the ground of this intuitive expectation. In fact, its possible truth and, even more, any suitable quantification of its validity would provide us with invaluable information about the approximation of Bayesian previsions by frequentistic ones which---as already explained---although cruder, are of easier evaluation. This circumstance comes to the fore, for example, within the so-called \emph{empirical Bayes} approach, which tries to justify partial replacement of orthodox Bayesian reasoning with frequentistic elements. See, e.g., Robbins (1956, 1964), Efron (2003) and Remark 2 in Section \ref{sect:merging} of the present paper.

The present work, which is part of a wide-ranging research, focuses on the discrepancy between posterior (predictive of $m$ future observations, respectively) distribution, given $n$ past observations, and the point mass at (the $m$-fold product of, respectively) the empirical distribution of the same past observations, when $n$ goes to infinity. The idea to compare a Bayesian inference to any of its frequentistic counterparts goes back, for different motives, to classical authors, such as
Laplace (1812), Poincar\'{e} (1912), Bernstein (1917), von Mises (1919, 1964), de Finetti (1929, 1930, 1937), Romanovsky (1931), and has had remarkable developments also in recent years, at least in two directions: The \emph{consistency} of Bayesian procedures from a frequentistic point of view, and the \emph{Bernstein-von Mises} phenomenon concerning a version of the central limit theorem for Bayesian estimators, in order to provide confidence regions connected with the aforesaid consistency issue. By way of example, see Schwartz (1965), Diaconis and Freedman (1986), Barron, Schervish and Wasserman (1999), Ghosal, Ghosh and van der Vaart (2000), Ghosh and Ramamoorthy (2003) for the consistency, and this last book together with Freedman (1999) for the Bernstein-von Mises phenomenon. To explain the connection with the present work, one must say that, especially in recent times, these fields of research have aligned themselves more and more with the interpretation of Bayesian inferences as procedures aimed at producing suitable estimators of unknown quantities, whose efficiency is checked from a frequentistic viewpoint. The product is to devoid both posterior and predictive distributions of their original meaning and role of conditional distributions, to reduce them to mere functions of the observations. Therefore, to appreciate the peculiarity of our work with respect to the aforesaid lines of research, one should thoroughly retrieve the Bayesian approach to statistical inference, in the spirit of the solution to the problem of \emph{inverse probabilities} provided by de Finetti (1929, 1930)
in his earliest papers on exchangeability. Nowadays, Doob (1949) is commonly credited as the author of the solution to a generalized form of the same problem, recalled in Theorem \ref{thm:Doob} of this paper. Indeed, if one reckons that the Bayesian way of thinking indicates, \emph{lato sensu}, the correct way of making statistical inference, it is fair to pursue the above-mentioned goals of approximating posterior and predictive distributions by more tractable laws---typically obtained by frequentistic procedures---depending only on past observations. Doob's theorem is then replaced by a statement concerning the almost sure (a.s.) convergence to zero of any weak probability distance (see Subsection \ref{sect:metric} below for more information) between the posterior distribution $q(\xisun, \cdot)$ and $\delta_{\empiric}$, the point mass at the empirical measure $\empiric := \frac{1}{n} \sum_{i=1}^n \delta_{\tilde{\xi}_i}$, as $n \rightarrow +\infty$. Successively, one can deduce the a.s.\! convergence to zero of any weak probability distance between the predictive distribution of $m$ future observations and the $m$-fold product $\empiric^m := \underbrace{\empiric \otimes \dots \otimes \empiric}_{\substack{m-\text{times}}}$, as $n \rightarrow +\infty$, for every $m \in \naturals$. Moreover, the main results in the present paper involve only \emph{finitary}---hence, empirically ascertainable---entities. In this respect, see Bassetti (2011) for the connection with finite exchangeable sequences. At this stage, one can appreciate the further step, made in Theorems \ref{thm:MainW}, \ref{thm:MainG} and \ref{thm:MainP}, to provide quantitative estimations of the error in the aforesaid approximations. More precisely, considering by way of example the comparison of $q(\xisun, \cdot)$ with $\delta_{\empiric}$, there are
a positive (non random) sequence $b_n$, going to infinity with $n$, and a suitable constant $L > 0$ such that, for every $\varepsilon, \eta > 0$, there exists some index $n_0 = n_0(\varepsilon, \eta) \in \naturals$ satisfying
$$
\rho\left(\left\{\max_{\nu \leq n \leq \nu+m} b_n \ud_{[\pms]}(q(\xisun), \delta_{\empiric}) \leq L + \varepsilon \right\}\right) \geq 1-\eta
$$
for every $\nu \geq n_0$ and $m \in \naturals$, where $\rho$ denotes the p.d. that makes the $\tilde{\xi}_n$'s exchangeable, and $\ud_{[\pms]}$ is a suitable probability distance to be specified in Subsection \ref{sect:metric}. See Remark \ref{rmk:finitistico} in Section \ref{sect:merging} for more explanation. Allied results, formulated in similar frameworks, can be found in Diaconis and Freedman (1990) and in Berti et. al. (2009).

Going back to posterior consistency, the analysis of the rapidity of convergence is usually expressed, like in Ghosal, Ghosh and van der Vaart (2000),  as
$$
\lim_{n \rightarrow +\infty} \pfrak_0^{\infty}\left(\left\{ q(\xisun, \{b_n \ud_{\pms}(\randommeasure, \pfrak_0) \geq M\}) \geq \varepsilon\right\}\right) = 0
$$
where $\ud_{\pms}$ stands for a suitable probability distance between probability laws on $\samplespace$ and $\pfrak_0^{\infty}$ denotes the p.d. of an infinite sequence of i.i.d. random element with common distribution $\pfrak_0$. This statement suffers a number of drawbacks substantially stemming from the co-existence of a p.d. with another thought of as ``true'', that are $\rho$ and $\pfrak_0^{\infty}$, respectively.

The last merit of our results is connected with the metric issue illustrated in Subsection \ref{sect:metric}. In fact, discrepancy between probability laws is here measured by means of probability metrics, which take account of any possible metric structure naturally attached to $\samplespace$, unlike the use of other measures of concentration (such as Kolmogorov-like ``distances'' and Kullback-Leibler divergences), which take maximum values, for example, in comparing point masses independently of any natural distance between the points of degeneracy. In this respect, it is worth mentioning significant works aimed at proving versions of the classical Glivenko-Cantelli theorem in the presence of probability distances to measure discrepancy between the empirical distribution and the ``true law''. See de Finetti (1933), Dudley (1969) and Yukic (1989).

As to the organization of the present work, the main results are formulated in Section \ref{sect:merging}, while Section \ref{sect:preliminary} recalls some preliminary facts about exchangeability, the Bayes-Laplace paradigm and the already mentioned metric issue. Finally, Section \ref{sect:proofs} contains the proofs of the new results.

\section{Preliminaries} \label{sect:preliminary}

A few preliminary notions, concepts and results are gathered in the present section to facilitate understanding of next developments. Subsection \ref{sect:exchangeable} describes the essentials of exchangeability with a view to its use in a general nonparametric framework. Subsection \ref{sect:BL} recalls the precise notions of prior, posterior and predictive distributions, and includes a statement on the limiting behavior of the last two as the number of observations goes to infinity. Finally, Subsection \ref{sect:metric} deals with some aspects about the metrization of both product spaces and spaces of probability measures (p.m.'s).

\subsection{Exchangeable observations} \label{sect:exchangeable}

Assume each observation takes values in a set $\samplespace$, a Borel subset of some Polish space $\hat{\samplespace}$,
and denote by $\pms$ \footnote{This notation is borrowed from de Finetti (1952).} the set of all p.m.'s on $(\samplespace, \mathscr{B}(\samplespace))$ where, as usual, for any topological space $\Space$, $\mathscr{B}(\Space)$ indicates the Borel class on $\Space$. Endow $\pms$ with the topology of weak convergence of p.m.'s and recall that, in view of the separability of $\samplespace$, $\mathscr{B}(\pms)$ is the same as the $\sigma$-algebra generated by the sets $\{\pfrak \in \pms\ |\ \pfrak(A) \in L\}$ as $A$ varies in $\mathscr{B}(\samplespace)$ and $L$ in $\mathscr{B}([0, 1])$. See, e.g., Proposition A2.5.IV in Daley and Vere-Jones (2003). Then, consider an infinite sequence of \emph{exchangeable} observations, in the sense that the probability distribution (p.d.) of each $k$-uple of distinct elements of the sequence depends only on $k$, for every $k \in \naturals$. With a view to next developments, identify this sequence of observations with the sequence $\xib = (\tilde{\xi}_1, \tilde{\xi}_2, \dots)$ of coordinate random elements of the usual topological product space $\samplespace^{\infty}$, endowed with its Borel $\sigma$-algebra $\mathscr{B}(\samplespace^{\infty})$. This way, the $j$-th component $\tilde{\xi}_j$ of $\xib$ is naturally
associated with the $j$-th observation, for any $j \in \naturals$. To complete the notation concerned with observable quantities, indicate by $\xisun$ and $\xisunm$ the vectors $(\tilde{\xi}_1, \dots, \tilde{\xi}_n)$ and $(\tilde{\xi}_{n+1}, \dots, \tilde{\xi}_{n+m})$, respectively, for every $n \in \naturals$ and $m \in \naturals \cup \{+\infty\}$. Now, it is well-known that under the above topological assumptions and ensuing measurability provisos, for any p.m. $\rho$ on $(\samplespace^{\infty}, \mathscr{B}(\samplespace^{\infty}))$ making the $\tilde{\xi}_j$'s exchangeable, there exists a unique p.m. $q$ on $(\pms, \mathscr{B}(\pms))$ such that \emph{de Finetti's representation}
\begin{equation} \label{eq:dF}
\rho(C) = \int_{\pms} \pfrak^{\infty}(C) q(\ud \pfrak)\ \ \ \ \ \ \ \ \ \ \ (C \in \mathscr{B}(\samplespace^{\infty}))
\end{equation}
holds true. See Aldous (1985) for a comprehensive treatment of exchangeability, included that
\begin{equation} \label{eq:dF2}
\rho(\{\empiric \Rightarrow \randommeasure,\ \text{as}\ n \rightarrow +\infty\}) = 1
\end{equation}
is in force, where $\Rightarrow$ denotes weak convergence of p.m.'s and $\randommeasure$ is a random p.m. on $(\samplespace, \mathscr{B}(\samplespace))$ having $q$ as p.d.. Recall that $\randommeasure$ is called a \emph{random p.m.} if it is a $\mathscr{B}(\samplespace^{\infty})/\mathscr{B}(\pms)$-measurable function from $\samplespace^{\infty}$ into $\pms$. Finally,
by the same de Finetti's theorem, $\randommeasure^{\infty}$ turns out to be a version of the regular conditional p.d. of $\xib$ given $\randommeasure$ or, equivalently, given either the tail $\sigma$-algebra $\mathscr{T}$ of $\xib$ or the $\sigma$-algebra $\mathscr{E}$ of the exchangeable events contained in $\mathscr{B}(\samplespace^{\infty})$.

\subsection{The Bayes-Laplace paradigm} \label{sect:BL}

The above statement regarding $\randommeasure^{\infty}$ as a version of conditional p.d. corresponds to the nonparametric form of the Bayes-Laplace paradigm for conditionally i.i.d. observations, when $q$ is the \emph{prior} p.d.. Bayesian statistical inferences on  $\randommeasure$ are based on a conditional p.d. for $\randommeasure$ given $\xisun$, the so-called \emph{posterior} p.d.. In view of the conditions stipulated in the previous subsection, well-known theorems (see, e.g., Theorems 6.3 and A1.2 in Kallenberg (2002)) can be applied to state the existence of a posterior p.d. given $\xisun$, say $q(\xisun) :=
q(\xisun, \cdot)$, that is:
\begin{enumerate}
\item[a)] $x \mapsto q(x, B)$ is $\sigma(\xisun)$-measurable, for every $B \in \mathscr{B}(\pms)$;
\item[b)] $B \mapsto q(x, B)$ is a p.m. on $(\pms, \mathscr{B}(\pms))$, for every $x \in \samplespace^n$;
\item[c)] $\int_C q(x, B) \rho^{(n)}(\ud x) = \rho(\{\xisun \in C, \randommeasure \in B\})$ holds for every $B \in \mathscr{B}(\pms)$ and $C \in \mathscr{B}(\samplespace^n)$, with $\rho^{(n)}(C) := \rho(C \times \samplespace^{\infty})$.
\end{enumerate}

According to the programme fixed in the introduction, prevision of future facts based on observed facts is at the core of the present paper. Such a kind of prevision relies on any conditional p.d. for $\xisunm$ given $\xisun$, generally named \emph{predictive}, a version of which, say $p_m(\xisun) := p_m(\xisun, \cdot)$, can be expressed by means of the posterior $q(\xisun)$ as
\begin{equation} \label{eq:predictive}
p_m(\xisun, C) = \int_{\pms} \pfrak^m(C) q(\xisun, \ud \pfrak)\ \ \ \ \ \ \ \ \ \ \ (C \in \mathscr{B}(\samplespace^m))
\end{equation}
for any $n \in \naturals$ and $m \in \naturals \cup \{+\infty\}$. Therefore, the role played by the posterior distribution may be important even if the main purpose is prevision of unknown facts based on observed facts.

Usually, center of attention of Bayesian statisticians is a random statistical model, seen as an unknown (or partially known) entity to be approached by means of statistical methods. In the above nonparametric setting, the role of model is played by
$\randommeasure^{\infty}$ and hence the following limit theorem may be of some importance. It can be derived from (\ref{eq:dF2}) in a rather direct way.
\begin{thm} \label{thm:Doob}
\emph{If} $\samplespace$ \emph{is a Borel subset of some Polish space and} $\rho$ \emph{is assessed like in} (\ref{eq:dF}), \emph{then}
\begin{eqnarray}
q(\xisun) \Rightarrow \delta_{\randommeasure}\ \ \ \text{as}\ n \rightarrow +\infty & \ \ \ \ \ (\rho-a.s.) \label{eq:DoobFin} \\
p_m(\xisun) \Rightarrow \randommeasure^m\ \ \ \text{as}\ n \rightarrow +\infty & \ \ \ \ \ (\rho-a.s.) \label{eq:DoobPred}
\end{eqnarray}
\emph{hold for every} $m \in \naturals \cup \{+\infty\}$.
\end{thm}

The proof can be found in Diaconis and Freedman (1986).

\subsection{The metric issue} \label{sect:metric}

The establishment of versions of (\ref{eq:DoobFin})-(\ref{eq:DoobPred}) in which $\delta_{\randommeasure}$ and $\randommeasure^m$ are replaced by $\delta_{\empiric}$ and $\empiric^m$, respectively, requires the introduction of suitable indices of discrepancy between p.m.'s, typically expressed as \emph{probability distances}. Actually, the problem of quantifying discrepancy between p.d.'s has to do with many areas of mathematics and has been dealt with from several viewpoints. See, for example, Rachev et al. (2013). This is why it is provided here a brief account of the distances that will be used in the rest of the work, motivating the choice by saying that the main aim is to situate the reasoning in a mathematical context as general as possible. In addition, the introduction of a distance is a prerequisite for the quantification of the phenomena of merging announced in Section \ref{sect:intro}. To get to the heart of the matter, notice that there are several cases in which the sample space $\samplespace$ is endowed with a natural distance $\ud_{\samplespace}$, as it happens, e.g., when $\samplespace = \reals^d$ and $\ud_{\samplespace}$ coincides with the Euclidean distance. On the contrary, it could happen that $\samplespace$ is specified only at the level of topological space, here supposed with the same features as at the beginning of Subsection \ref{sect:exchangeable}. In that case, recall that the topology $\mathcal{T}_{\samplespace}$ can be thought of as generated by a totally bounded metric $\ud_{\samplespace}^{'}$ on $\samplespace$, that in the specific case plays the role of mere mathematical device. A proof of the existence of $\ud_{\samplespace}^{'}$ is contained in Sections 2.6 and 2.8 of Dudley (2002). Situations of this kind, in which it is immaterial whether one adopts one metric or another, are rather common in the statistical analysis of real problems involving, for instance, qualitative characters or infinite-dimensional mathematical objects belonging to spaces whose topologies are characterized through collections of neighborhoods. All these circumstances come into play at the same time in which $[\samplespace]$, $[[\samplespace]]$, $\samplespace^m$ and  $[\samplespace^m]$ require a metric structure, so that the remaining part of this subsection will be devoted to an illustration of this aspect.

When a distinguished metric $\ud_{\samplespace}$ is given, there are many significant distances on $[\samplespace]$ whose definition rests crucially on $\ud_{\samplespace}$, but the present work makes use only of the \emph{Prokhorov} and the \emph{Gini-Monge-Wasserstein} metric. The former is defined by
$$
\ud_{[\samplespace]}^{(P)}(\mu_1, \mu_2) := \inf\{\epsilon > 0\ | \ \mu_1(B) \leq \mu_2(B^{\epsilon}) + \epsilon,\ \forall \ B \in \mathscr{B}(\samplespace)\}
$$
where $B^{\epsilon} := \{x \in \samplespace\ |\ \ud_{\samplespace}(x, y) < \epsilon\ \text{for\ some}\ y \in B\}$. As to the latter, given $\mu_1, \mu_2 \in [\samplespace]$, let $\mathcal{F}(\mu_1, \mu_2)$ stand for the class of all p.m.'s on $(\samplespace^2, \mathscr{B}(\samplespace^2))$ with $i$-th marginal equal to $\mu_i$, $i = 1, 2$. If, for some $p \in [1, +\infty)$, $\mu_i \in [\samplespace]_p := \big{\{}\mu \in [\samplespace]\ \big{|}\ \int_{\samplespace} [\ud_{\samplespace}(x, x_0)]^p \mu(\ud x) < +\infty\ \text{for\ some}\ x_0 \in \samplespace\big{\}}$ for $i = 1, 2$, the \emph{Gini-Monge-Wasserstein} distance of order $p$ between $\mu_1$ and $\mu_2$ is defined to be
$$
\ud_{[\samplespace]}^{(G_p)}(\mu_1, \mu_2) := \inf_{\gamma \in \mathcal{F}(\mu_1, \mu_2)}\left(
\int_{\samplespace^2} [\ud_{\samplespace}(x, y)]^p\ \gamma(\ud x \ud y)\right)^{1/p}\ .
$$
The definition of these two metrics can be extended (with respective notations $\ud_{[[\samplespace]]}^{(P)}$ and $\ud_{[[\samplespace]]}^{(G_p)}$) to the space $[[\samplespace]]$ by replacing, in the last two formulas, $\samplespace$ with $[\samplespace]$ and by making the proviso that, in the expression of $\ud_{[[\samplespace]]}^{(P)}$ ($\ud_{[[\samplespace]]}^{(G_p)}$, respectively), $\ud_{\samplespace}$ is replaced with $\ud_{[\samplespace]}^{(P)}$ ($\ud_{[\samplespace]}^{(G_p)}$, respectively). Apropos of $\samplespace^m$, we observe that any usual product metric, such as $\left(\sum_{i=1}^m [\ud(x_i, y_i)]^p \right)^{1/p}$ with $p \in [1, +\infty)$, does not match with the assumption of exchangeability, due to the lack of invariance under permutation of the coordinates of each single vector. Therefore, inspired by the original works by Gini (1914) and Leti (1961, 1962), we here propose to replace $\samplespace^m$ with its quotient $\samplespace^m_{\sigma} := \samplespace^m/\sim$, where $\sim$ stands for the equivalence relation that identifies any vector $(x_1, \dots, x_m)$ with $(x_{\tau(1)}, \dots, x_{\tau(m)})$ for every $m$-permutation $\tau$, and to metrize $\samplespace^m_{\sigma}$ in such a way that the mapping $\samplespace^m_{\sigma} \ni [(x_1, \dots, x_m)] \mapsto  \frac{1}{m}\sum_{i=1}^m \delta_{x_i} \in [\samplespace]$ turns out to be an isometry. This plan can be carried out consistently with the metrization of $[\samplespace]$ described above, leading to
$$
\ud_{\samplespace^m_{\sigma}}^{(\star)}\left([(x_1, \dots, x_m)], [(y_1, \dots, y_m)]\right) := \ud_{[\samplespace]}^{(\star)}\left(\frac{1}{m}\sum_{i=1}^m \delta_{x_i}, \frac{1}{m}\sum_{i=1}^m \delta_{y_i}\right)
$$
where $[(x_1, \dots, x_m)]$ denotes the equivalence class of $(x_1, \dots, x_m)$ and $\star$ stands either for $P$ or $G_p$. Finally, the metrization of $[\samplespace^m_{\sigma}]$ parallels that of $[[\samplespace]]$, with the proviso that in the expression of $\ud_{[\samplespace^m_{\sigma}]}^{(P)}$ ($\ud_{[\samplespace^m_{\sigma}]}^{(G_p)}$, respectively),
$\ud_{\samplespace^m_{\sigma}}^{(P)}$ ($\ud_{\samplespace^m_{\sigma}}^{(G_p)}$, respectively) appears in the place of
$\ud_{\samplespace}$.

The last part of this subsection re-examines the previous picture in the event that $\samplespace$ is simply given as Borel subset of some Polish space, and the adoption of one metric or another is considered as immaterial. Upon imposing a totally bounded metric $\ud^{'}_{\samplespace}$ on $\samplespace$ as above,
recall that there exists a countable collection of $\ud^{'}_{\samplespace}$-uniformly continuous, $[0, 1]$-valued functions
forming a determining class for weak convergence in $\pms$. See, e.g., Theorem 9.1.5 in Stroock (2011). Moreover, observe that this class can be chosen, without any loss of generality, with the following additional properties:
i) every function is $\ud^{'}_{\samplespace}$-Lipschitz continuous; ii) the entire class is dense, with respect to the sup norm, in the space of $\ud^{'}_{\samplespace}$-uniformly continuous, $[0, 1]$-valued functions. See Subsection \ref{sect:detclass} for a proof of this fact. Thus, upon denoting by $\{g_k\}_{k \geq 1}$ such a class and by $\|g_k\|_{BL}$ the norm $\|g_k\|_{\infty} + \|g_k\|_{\text{Lip}} := \sup_{x \in \samplespace} |g_k(x)| + \sup_{x \neq y} [|g_k(x) - g_k(y)|
/\ud_{\samplespace}^{'}(x, y)]$, define
$$
\ud_{\pms}^{(W)}(\mu_1, \mu_2) := \sum_{k=1}^{\infty} \frac{1}{2^k } \Big{|} \int_{\samplespace} g_k^{\ast}(x) \mu_1(\ud x) - \int_{\samplespace} g_k^{\ast}(x) \mu_2(\ud x) \Big{|}
$$
to be the desired metric on $\pms$, with $g_k^{\ast}(x) := g_k(x)/\|g_k\|_{BL}$. To metrize $[[\samplespace]]$, observe that $(\pms, \ud_{\pms}^{(W)})$ is separable, and repeat step by step the above construction with a new sequence $\{h_k\}_{k \geq 1}$ of $\ud_{\pms}^{(W)}$-Lipschitz, $[0, 1]$-valued functions. Finally, after noting that the countable family $\{\prod_{i=1}^m g_{k_i}(x_i)\}_{k_1, \dots, k_m \in \naturals}$ forms a determining class for weak convergence in $[\samplespace^m]$ (see again Subsection \ref{sect:detclass}), use the metric
$$
\ud_{[\samplespace^m]}^{(W)}(\nu_1, \nu_2) := \sum_{k_1, \dots, k_m \in \naturals} \frac{1}{2^{k_1 + \dots + k_m}} \Big{|} \int_{\samplespace^m} \left[\prod_{i=1}^m g_{k_i}^{\ast}(x_i)\right] \nu_1(\ud x) - \int_{\samplespace^m} \left[\prod_{i=1}^m g_{k_i}^{\ast}(x_i)\right] \nu_2(\ud x)\Big{|}
$$
to compare predictive distributions.

\section{Main results} \label{sect:merging}

The new results of this paper will be presented in the form of four statements. The first one, because of its qualitative nature, plays an introductory role by providing an analogous version of (\ref{eq:DoobFin})-(\ref{eq:DoobPred}), with $\empiric$ in the place of   $\randommeasure$. To achieve this aim, it will be necessary to consider two generic metrizations of weak convergence of p.m.'s on $\pms$ and $\samplespace^m$, denoted by $\ud_{[[\samplespace]]}$ and $\ud_{[\samplespace^m]}$, respectively.
\begin{thm} \label{thm:DoobMetrico}
\emph{If} $\samplespace$ \emph{is a Borel subset of some Polish space and} $\rho$ \emph{is given like in} (\ref{eq:dF}), \emph{then}
\begin{eqnarray}
\ud_{[[\samplespace]]}\big{(} q(\xisun), \delta_{\empiric} \big{)} \rightarrow 0\ \ \ \text{as}\ n \rightarrow +\infty & \ \ \ \ \ (\rho-a.s.) \label{eq:DoobQuantFin} \\
\ud_{[\samplespace^m]}\big{(} p_m(\xisun), \empiric^m\big{)} \rightarrow 0\ \ \ \text{as}\ n \rightarrow +\infty & \ \ \ \ \ (\rho-a.s.) \label{eq:DoobQuantPred}
\end{eqnarray}
\emph{hold true for every} $m \in \naturals \cup \{+\infty\}$.
\end{thm}
A proof is contained in Subsection \ref{sect:ProofDoobMetrico}. The next three theorems improve the last one by providing rates of approach to zero of the distances obtained by replacing the generic ones in (\ref{eq:DoobQuantFin})-(\ref{eq:DoobQuantPred}) with specific definitions given in Subsection \ref{sect:metric}. Indeed, the expression of any rate will be influenced not only by the probabilistic framework, encapsulated in the p.d. $\rho$, but also by the specific metric structure attached to the spaces $\samplespace$, $\samplespace^m$, $\pms$ and $[\pms]$. The value of the following results rests, above all, on the fact that the rates of approach to zero are \emph{deterministic}, and hence known to the statistician before getting the data.
\begin{thm} \label{thm:MainW}
\emph{Assume that} $\samplespace$, \emph{a Borel subset of some Polish space, is metrized by a totally bounded distance} $\ud_{\samplespace}^{'}$, \emph{and that the spaces} $[\pms]$ \emph{and} $[\samplespace^m]$ \emph{are endowed with the metrics} $\ud_{[\pms]}^{(W)}$ \emph{and} $\ud_{[\samplespace^m]}^{(W)}$ \emph{respectively, as in Subsection} \ref{sect:metric}. \emph{Then, if} $\rho$ \emph{is given like in} (\ref{eq:dF}),
\begin{eqnarray}
\limsup_{n \rightarrow \infty} \sqrt{\frac{n}{\log\log n}} \ud_{[[\samplespace]]}^{(W)}\left(q(\xisun), \delta_{\empiric} \right) &\leq& \sqrt{2}\ \ \ \ \ \ \ \ \ (\rho-a.s.) \label{eq:MainDoobFin} \\
\limsup_{n \rightarrow \infty} \sqrt{\frac{n}{\log\log n}} \ud_{[\samplespace^m]}^{(W)}\left(p_m(\xisun), \empiric^m \right) &\leq& \sqrt{2} m \ \ \ \ \ \ (\rho-a.s.) \label{eq:MainDoobPred}
\end{eqnarray}
\emph{hold for every} $m \in \naturals$.
\end{thm}
The proof is developed in Subsection \ref{sect:ProofMainDoob}. Finally, the last three statements deal with the case in which the topology $\mathcal{T}_{\samplespace}$ is given in terms of some natural distance $\ud_{\samplespace}$ on $\samplespace$, starting from the noteworthy case $\samplespace = \reals$ and $\ud_{\samplespace}(x_1, x_2) := |x_1 - x_2|$ (Euclidean distance). Therefore, when the metric framework described in Subsection \ref{sect:metric} is based on the Gini-Monge-Wasserstein distance of order 1, one has the following
\begin{thm} \label{thm:MainG}
\emph{After choosing} $\ud_{[\reals]}^{(G_1)}$ \emph{as distance on} $[\reals]$, \emph{metrize} $[[\reals]]$, $\reals^m_{\sigma}$ \emph{and} $[\reals^m_{\sigma}]$ \emph{with} $\ud_{[[\reals]]}^{(G_1)}$, $\ud_{\reals^m_{\sigma}}^{(G_1)}$ \emph{and} $\ud_{[\reals^m_{\sigma}]}^{(G_1)}$ \emph{respectively, as in Subsection} \ref{sect:metric}. \emph{Moreover, given} $\rho$ \emph{as in} (\ref{eq:dF}), \emph{assume that} $\inteone |x|^{2+\delta} \overline{\pfrak}(\ud x) < +\infty$ \emph{obtains for some} $\delta > 0$, \emph{where} $\overline{\pfrak}(B) := \int_{[\reals]} \pfrak(B) q(\ud \pfrak)$ \emph{for every} $B \in \mathscr{B}(\reals)$. \emph{Then, putting} $\tilde{\textsf{F}}(x) := \randommeasure((-\infty, x])$, \emph{one has}
\begin{eqnarray}
\limsup_{n \rightarrow \infty} \sqrt{\frac{n}{\log\log n}} \ud_{[[\samplespace]]}^{(G_1)}\left(q(\xisun), \delta_{\empiric} \right) &\leq& \int_{\reals} \sqrt{2 \tilde{\textsf{F}}(x)[1 - \tilde{\textsf{F}}(x)]} \ud x \ \ \ (\rho-a.s.)  \label{eq:MainGiniFin} \\
\limsup_{n \rightarrow \infty} \sqrt{\frac{n}{\log\log n}} \ud_{[\samplespace^m_{\sigma}]}^{(G_1)}\left(p_m(\xisun), \empiric^m \right) &\leq& \int_{\reals} \sqrt{2 \tilde{\textsf{F}}(x)[1 - \tilde{\textsf{F}}(x)]} \ud x \ \ \ (\rho-a.s.)
\label{eq:MainGiniPred}
\end{eqnarray}
\emph{for every} $m \in \naturals$.
\end{thm}
For the proof see Subsection \ref{sect:ProofMainGini}. Here, it is worth mentioning the more convenient bound
$$
\int_{\reals} \sqrt{2 \tilde{\textsf{F}}(x)[1 - \tilde{\textsf{F}}(x)]} \ud x \leq \left(8\int_{\reals} \left[ |x| + \frac{|x|^{2+\epsilon}}{(2+\epsilon)} \right] \randommeasure(\ud x) \right)^{1/2}
$$ 
valid for every $\epsilon \in (0, \delta]$, that is displayed in the proof itself.

The next statement is concerned with the Prokhorov distance for probabilities on an abstract space which, in the present case,
satisfies the slightly more restrictive condition of being Borel subset of a \emph{locally compact} Polish space $\hat{\samplespace}$. Its thesis turns out to be less sharp than the previous ones for various reasons partially discussed in Remark \ref{rmk:prokhorov} at the end of this section. Here, suffice it to mention that the main hypothesis depends crucially on a discretization of $\randommeasure$, based on the fact that, thanks to the topological characterization of $\hat{\samplespace}$, one can deduce from Theorem 2.8.1 and Problem 2.8.6 in Dudley (2002) the existence of an increasing sequence of compact subsets $\hat{\mathbb{K}}_m$ converging to $\hat{\samplespace}$, and that each compact $\hat{\mathbb{K}}_m$ admits a partition $\{A_{m, j}\}_{j = 1, \dots, k_m}$ with $\text{diam}(A_{m, j}) \leq 1/m$. Then, putting $A_{m, k_m + 1} := \hat{\mathbb{K}}_m^c$ for all $m \in \naturals$ and considering
$$
\Pi_r(\pfrak) := \liminf_{m \rightarrow +\infty} \sum_{j=1}^{k_m + 1} [\pfrak(A_{m, j} \cap \samplespace) (1 - \pfrak(A_{m, j} \cap \samplespace))]^{1/r} \ \ \ \ \ \ \ (\pfrak \in \pms)
$$
for $r > 2$ pave the way for the formulation of
\begin{thm} \label{thm:MainP}
\emph{Assume that} $\samplespace$ \emph{is a Borel subset of a locally compact Polish space and metrize} $\pms$, $[\pms]$, $\samplespace_{\sigma}^m$ \emph{and} $[\samplespace_{\sigma}^m]$ \emph{with} $\ud_{\pms}^{(P)}$, $\ud_{[\pms]}^{(P)}$, $\ud_{\samplespace_{\sigma}^m}^{(P)}$  \emph{and} $\ud_{[\samplespace_{\sigma}^m]}^{(P)}$, \emph{respectively. Moreover, given}
$\rho$ \emph{as in} (\ref{eq:dF}), \emph{suppose that} $\Pi_r(\randommeasure) \in \mathrm{L}^1(\samplespace^{\infty}, \mathscr{B}(\samplespace^{\infty}), \rho)$ \emph{for some} $r > 2$. \emph{Then},
\begin{eqnarray}
\limsup_{n \rightarrow \infty} \left(\frac{n}{\log n}\right)^{1/4} \ud_{[[\samplespace]]}^{(P)}\left(q(\xisun), \delta_{\empiric} \right) &\leq& Y(\randommeasure) \ \ \ \ \ \ \ \ \ \ \ \ (\rho-a.s.) \label{eq:MainProkFin} \\
\limsup_{n \rightarrow \infty} \left(\frac{n}{\log n}\right)^{1/8} \ud_{[\samplespace^m_{\sigma}]}^{(P)}\left(p_m(\xisun), \empiric^m \right) &\leq& \sqrt{\frac{3}{2} Y(\randommeasure)} \ \ \ \ \ \ (\rho-a.s.) \label{eq:MainProkPred}
\end{eqnarray}
\emph{are valid for every} $m \in \naturals$, \emph{where} $Y(\randommeasure) := \left(\frac{3}{2} \limsup_{n \rightarrow +\infty} \sqrt{\frac{n}{\log n}} \ud_{\pms}^{(P)}(\randommeasure, \empiric)\right)^{1/2}$ \emph{is finite} $\rho-a.s.$.
\end{thm}
The proof is deferred to Subsection \ref{sect:ProofMainProk}, whilst the discussion here focuses on a new upper bound for $\Pi_r(\pfrak)$, valid when $\hat{\samplespace} = \rd$ and $\pfrak$ belongs to the distinguished class of probabilities
\begin{eqnarray}
\pms_{\ast} := \Big{\{} \pfrak(A) &=& \lambda \sum_{i = 1}^N p_i \delta_{x_i}(A) + (1 - \lambda) \int_A f(x) \ud x \ \text{for\ all}\ A \in \mathscr{B}(\samplespace), \nonumber \\
&& \text{for\ some}\ \lambda \in [0, 1],\ N \in \naturals,\ p_1, \dots, p_N \in [0, 1]\ \text{with}\ \sum_{i=1}^N p_i = 1 \nonumber \\
&& x_1, \dots, x_N \in \samplespace,\ f : \samplespace \rightarrow [0, \infty)\ \text{with}\ \int_{\samplespace} f(x)\ud x =1 \Big{\}}\ . \nonumber
\end{eqnarray}
\begin{prop} \label{prop:PropP}
\emph{Let} $\samplespace$ \emph{be either} $\rd$ \emph{or a regular} \footnote{This class contains the most common examples used in practical statistics, such as semi-spaces, cubes, balls, and so on. Typically, all bounded domains with Lipschitz boundary are included. The most general assumptions on $\samplespace$ are the so-called \emph{cone conditions} according to 4.6 in Adams and Fournier (2003).} \emph{sub-domain of} $\rd$. \emph{If} $q(\pms_{\ast}) = 1$ \emph{and both} $[\tilde{\lambda}]^{1/r} \sum_{i=1}^{\tilde{N}} [p_i]^{1/r}$ \emph{and} $[1 - \tilde{\lambda}]^{1/r} \sum_{|\boldsymbol{\alpha}| = 0}^{\frac{r-1}{r}d} \int_{\samplespace} |D^{\boldsymbol{\alpha}} (\tilde{f}(x))^{1/r}| \ud x$ \emph{belong to} $\mathrm{L}^1(\samplespace^{\infty}, \mathscr{B}(\samplespace^{\infty}), \rho)$ \emph{for some} $r > 2$ \emph{for which} $\frac{r-1}{r}d \in \naturals$, \emph{then
also} $\Pi_r(\randommeasure)$ \emph{belongs to} $\mathrm{L}^1(\samplespace^{\infty}, \mathscr{B}(\samplespace^{\infty}), \rho)$.
\end{prop}
See Subsection \ref{sect:sobolev} for a proof of this proposition.

The presentation of the main results terminates with some relevant comments.
\begin{remark}
It is worth noticing that, in the noteworthy case of a finite $\samplespace$, all the distances on $\pms$ mentioned in Subsection \ref{sect:metric} turn out to be metrically equivalent. Therefore, when $\ud_{[\samplespace^m]}^{(W)}$ in (\ref{eq:MainDoobPred}) is replaced by either $\ud_{[\samplespace^m]}^{(P)}$ or $\ud_{[\samplespace^m]}^{(G_p)}$, the bound therein remains valid up to a multiplicative factor, generally depending on the cardinality of $\samplespace$. As to (\ref{eq:MainDoobFin}), the replacement of $\ud_{[[\samplespace]]}^{(W)}$ with either  $\ud_{[[\samplespace]]}^{(P)}$ or $\ud_{[[\samplespace]]}^{(G_p)}$ is feasible in view of the following facts: First, $\ud_{[[\samplespace]]}^{(P)} \leq \sqrt{\frac{3}{2} \ud_{[[\samplespace]]}^{(G_p)}}$ for all $p \geq 1$ (see Problem 5, Section 11.3 of Dudley (2002)). Second, a well-known Kantorovich-Rubinstein representation of $\ud_{[[\samplespace]]}^{(G_1)}$ can be carried out exactly as in Subsection  \ref{sect:ProofMainDoob}.
\end{remark}

\begin{remark} \label{rmk:finitistico}
Here is a formal argument to retrieve the more practical reformulation of the previous results, mentioned in the Introduction. One confines oneself to considering the bounds (\ref{eq:MainDoobFin}),(\ref{eq:MainGiniFin}) and (\ref{eq:MainProkFin}) concerning the posterior distribution, which can be reduced to the form
$$
\limsup_{n \rightarrow \infty} b_n \ud_{[[\samplespace]]}\left(q(\xisun), \delta_{\empiric} \right) \leq \tilde{Y}
$$
for some random number $\tilde{Y} = Y(\randommeasure)$. To carry out the desired reformulation, one first fixes $\eta > 0$ and determines an $L > 0$ such that $\rho(\tilde{Y} > L) \leq \eta/2$ and, as a consequence,
$$
1 - \eta/2 \leq \rho\left(\left\{\limsup_{n \rightarrow \infty} b_n \ud_{[[\samplespace]]}\left(q(\xisun), \delta_{\empiric} \right) \leq L\right\}\right)\ .
$$
At this stage, for every $\varepsilon > 0$, one can determine $n_0 = n_0(\varepsilon, \eta)$ such that
$$
1 - \eta \leq \rho\left(\left\{\max_{n \geq \nu} b_n \ud_{[[\samplespace]]}\left(q(\xisun), \delta_{\empiric} \right) \leq L
+ \varepsilon \right\}\right)
$$
for every $\nu \geq n_0$. To conclude, suffice it to notice that the right-hand side is not greater than \\ $\rho\left(\left\{\max_{\nu \leq n \leq \nu + m} b_n \ud_{[[\samplespace]]}\left(q(\xisun), \delta_{\empiric} \right) \leq L
+ \varepsilon \right\}\right)$ for every $m \in  \naturals$.
\end{remark}

\begin{remark} \label{rmk:empBayes}
Here is an application to the evaluation of the error provoked by the plugging of frequentistic components in a Bayesian inference, following the empirical Bayes approach. To clear the field from unessential complications, one considers the basic problem of estimating the expectation $m_g(\randommeasure) := \int_{\samplespace} g(x) \randommeasure(\ud x)$ on the basis of $\xisun$, when $g$ is a fixed measurable function from $\samplespace$ to $\reals$ such that $\int_{\samplespace} |g(x)|^{2 + \delta} \overline{\pfrak}(\ud x) < +\infty$ for some $\delta > 0$, with $\overline{\pfrak}(B) := \int_{[\reals]} \pfrak(B) q(\ud \pfrak)$. Assume that the Bayesian estimator is $\tilde{B}_n := \int_{\pms} m_g(\pfrak) q(\xisun, \ud \pfrak)$, but that the statistician interested in its evaluation is unable to specify any prior distribution, so that she/he decides to fall back on the most obvious frequentistic solution $\tilde{g}_n := \int_{\samplespace} g(x) \empiric(\ud x)$. A question could arise as to the proximity of this convenient arrangement to the Bayesian solution $\tilde{B}_n$. One naturally expects the answer to depend on the prior $q$ or, at least, on some particular aspect of it. As an asymptotic measure of proximity here one chooses to analyze the behavior of $b_n |\tilde{B}_n - \tilde{g}_n|$ as $n \rightarrow +\infty$, for some suitable sequence $\{b_n\}_{n \geq 1}$ going to infinity with $n$. To this end, notice that
\begin{eqnarray}
\tilde{B}_n &=& \ee\left[\ee\left(g(\tilde{\xi}_{n+1})\ |\ \randommeasure\right)\ \big{|}\ \xisun\right]
= \ee\left[\ee\left(g(\tilde{\xi}_{n+1})\ |\ \randommeasure, \xisun\right)\ \big{|}\ \xisun\right] \nonumber\\
&=& \ee\left(g(\tilde{\xi}_{n+1})\ |\ \xisun\right) = \int_{\reals} x p_1(\xisun) \circ g^{-1}(\ud x) \nonumber
\end{eqnarray}
and, from the definition of the Gini-Monge-Wasserstein distance, one gets
$$
|\tilde{B}_n - \tilde{g}_n| \leq \ud_{[\reals]}^{(G_1)}\big{(} p_1(\xisun) \circ g^{-1}, \empiric \circ g^{-1} \big{)}\ .
$$
Whence, a straightforward application of (\ref{eq:MainGiniPred}) yields
$$
\rho\left(\left\{ \limsup_{n \rightarrow +\infty} \sqrt{\frac{n}{\log\log n}} |\tilde{B}_n - \tilde{g}_n| \leq \int_{\reals} \sqrt{2 \tilde{\textsf{F}}_g(x)[1 - \tilde{\textsf{F}}_g(x)]} \ud x \right\}\right) = 1
$$
with $\tilde{\textsf{F}}_g(x) := \randommeasure \circ g^{-1}((-\infty, x])$. For a more practical interpretation of this result, one can go back to the previous remark.
\end{remark}

\begin{remark} \label{rmk:prokhorov}
The difficulties connected with Theorem \ref{thm:MainP} and Proposition \ref{prop:PropP} are due to the generality of the space $\samplespace$. In fact, the bounds exhibited in (\ref{eq:MainProkFin})-(\ref{eq:MainProkPred}), being of an implicit nature, are not so useful. Nonetheless, they have the value to connect the rapidity of merging of Bayesian inferences with their frequentistic counterparts to the speed of Glivenko-Cantelli convergence. Indications about that phenomenon can be found in Dudley (1969) and Yukich (1989), although a complete characterization is still lacking even when $\samplespace \subset \rd$. Further improvements could be obtained by a more precise investigation on the exit times of a simple random walk with respect to suitable curvilinear boundaries, which will constitute the subject of a forthcoming paper. Finally, Proposition \ref{prop:PropP} presents some sufficient conditions which can be of simpler verification, for example, in connection with distinguished priors, whose support is included in the space of probability density functions on $\samplespace$. Noteworthy examples are the so-called ``models for density estimation'', which include, for example, the \emph{mixture models} proposed by Lo (1984).
\end{remark}

\section{Proofs} \label{sect:proofs}

Gathered in this section are the proofs of the statements formulated in the previous ones. The present one is split into a few subsections, the first of which completes the arguments used in Subsection \ref{sect:metric} to characterize the classes of test functions for the metrics $\ud_{[[\samplespace]]}^{(W)}$ and $\ud_{[\samplespace^m]}^{(W)}$. The second provides the proof of Theorem \ref{sect:ProofDoobMetrico}. The third contains a complement to Theorem 4 in Section 10.3 of Chow and Teicher (1997) on the maximum of normed sums. Finally, in the last subsections, one proves Theorems \ref{thm:MainW}-\ref{thm:MainP}
and Proposition \ref{prop:PropP}.

\subsection{On certain weak convergence-determining classes} \label{sect:detclass}

Define $\mathcal{U}_b^{\ud_{\samplespace}^{'}}(\samplespace; [0, 1])$ to be the class of $\ud_{\samplespace}^{'}$-uniformly continuous functions, with values in $[0, 1]$. According to Lemma 9.1.4 in Stroock (1999), there is a countable subclass $\{u_k\}_{k \geq 1}$ which is dense in $\mathcal{U}_b^{\ud_{\samplespace}^{'}}(\samplespace; [0, 1])$ with respect to the $\sup$-norm. Since the completion $\hat{\samplespace}$ of $\samplespace$ is compact (c.f. Theorem 2.8.2 in Dudley (2002)), each $u_k$ can be extended to a uniformly continuous function $\hat{u}_k$ on $\hat{\samplespace}$. Then, from Theorem 11.2.4 in Dudley (2002) on density of Lipschitz functions, each $\hat{u}_k$ can be uniformly approximated by a suitable sequence $\{\hat{u}_{k,n}\}_{n \geq 1}$ of $\ud_{\hat{\samplespace}}^{'}$-Lipschitz continuous, $[0, 1]$-valued functions. Now, define $u_{k,n}$ to be the restriction to $\samplespace$ of $\hat{u}_{k,n}$, for every $k, n \in \naturals$. At this stage, the desired class $\{g_k\}_{k \geq 1}$ can be obtained by re-enumerating $\{u_{k, n}\}_{k,n \in \naturals}$.

To prove that $\{\prod_{i=1}^m g_{k_i}(x_i)\}_{k_1, \dots, k_m \in \naturals}$ represents a determining class for weak convergence of p.m.'s on $\samplespace^m$, first notice that $\nu_n^{(m)} \Rightarrow \nu^{(m)}$ as $n \rightarrow +\infty$, for $\nu_n^{(m)}, \nu^{(m)} \in [\samplespace^m]$, is equivalent to 
$$
\int_{\samplespace^m} \left[\prod_{i=1}^m f_i(x_i)\right] \nu_n^{(m)} (\ud x_1 \dots \ud x_m) \rightarrow \int_{\samplespace^m} \left[\prod_{i=1}^m f_i(x_i)\right] \nu^{(m)} (\ud x_1 \dots \ud x_m)
$$ 
as $n \rightarrow +\infty$, for every $f_1, \dots, f_m \in \mathcal{U}_b^{\ud_{\samplespace}^{'}}(\samplespace; [0, 1])$. See, e.g., Corollary 1.4.5 in van der Vaart and Wellner (1996). To complete the argument, combine the density of $\{g_k\}_{k \geq 1}$ in $\mathcal{U}_b^{\ud_{\samplespace}^{'}}(\samplespace; [0, 1])$ with Lemma 1 of Section 27 in Billingsley (1995) concerning the difference of products of complex numbers.

\subsection{Proof of Theorem \ref{thm:DoobMetrico}} \label{sect:ProofDoobMetrico}

Apply the triangle inequality to get
\begin{eqnarray}
\ud_{[[\samplespace]]}\big{(} q(\xisun), \delta_{\empiric} \big{)} &\leq& \ud_{[[\samplespace]]}\big{(} q(\xisun), \delta_{\randommeasure} \big{)} + \ud_{[[\samplespace]]}\big{(} \delta_{\randommeasure}, \delta_{\empiric} \big{)} \nonumber \\
\ud_{[\samplespace^m]}\big{(} p_m(\xisun), \empiric^m\big{)} &\leq& \ud_{[\samplespace^m]}\big{(} p_m(\xisun), \randommeasure^m\big{)} + \ud_{[\samplespace^m]}\big{(} \randommeasure^m, \empiric^m\big{)} \nonumber
\end{eqnarray}
where, by virtue of Theorem \ref{thm:Doob}, both terms $\ud_{[[\samplespace]]}\big{(} q(\xisun), \delta_{\randommeasure} \big{)}$ and $\ud_{[\samplespace^m]}\big{(} p_m(\xisun), \randommeasure^m\big{)}$ go to zero with probability one, as $n$ goes to infinity and for every fixed $m \in \naturals \cup \{+\infty\}$. To deduce the validity of (\ref{eq:DoobQuantFin}),
combine Theorem 11.3.3 in Dudley (2002) on equivalence of probability metrics with the obvious identity $\ud_{[[\samplespace]]}^{(P)}\big{(} \delta_{\pfrak_1}, \delta_{\pfrak_2} \big{)} = \min\left\{1, \ud_{\pms}(\pfrak_1, \pfrak_2)\right\}$ (c.f. Section 11.3 of the same book) to prove that (\ref{eq:dF2}) entails $\delta_{\empiric} \Rightarrow \delta_{\randommeasure}$ with probability one, as $n$ goes to infinity. Finally, to deduce (\ref{eq:DoobQuantPred}), invoke (\ref{eq:dF2}) once again and apply in combination Theorem 4.29 in Kallenberg (2002) and Theorem 2.8(ii) of Billingsley (1999) on weak convergence in product spaces.

\subsection{Complement to Teicher's theorem on the maximum of normed sums} \label{sect:teicher}

In each of the proofs of Theorems \ref{thm:MainW}-\ref{thm:MainP} there is a step based on the following theorem by Teicher. See
Theorem 4 in Section 10.3 of Chow and Teicher (1997). This result is concerned with a sequence $\{X_n, n \geq 1\}$ of i.i.d. random numbers, defined on some probability space $\probabilityspace$: If $\ee[X_1] = 0$ and $\ee[|X|^r] < +\infty$ for some $r > 2$, where $\ee$ denotes expectation, then
\begin{equation} \label{eq:teicherpiu2}
\ee\left[\sup_{n > e^e} \frac{\big{|}\sum_{i}^n X_i\big{|}^r}{(n \log\log n)^{r/2}}\right] \leq \alpha_0(r) \sigma^r + \alpha_1(r) (\sigma^r)^{1- \lceil r \rceil} \left(\ee[|X|^r]\right)^{\lceil r \rceil}
\end{equation}
obtains, where $\alpha_0(r), \alpha_1(r)$ are suitable constants which do not depend on $\pp$, $\sigma^2 := \ee[X^2]$ and $\lceil r \rceil$ stands for $\inf\{n \in \naturals\ |\ n \geq r\}$. The bound specified in (\ref{eq:teicherpiu2}) is obtained simply by detailing the original Teicher proof.

In this very same setting, one can deduce another useful bound by simply combining the proof given by Teicher with the so-called Rosenthal inequality (see, e.g., Section 2.3 in Petrov (1995)), namely
\begin{equation} \label{eq:newteicherpiu2}
\ee\left[\sup_{n > e} \frac{\big{|}\sum_{i}^n X_i\big{|}^r}{(n \log n)^{r/2}}\right] \leq \beta_0(r) \sigma^r + \beta_1(r)\ee[|X|^r]
\end{equation}
where again $\beta_0(r), \beta_1(r)$ are suitable constants independent of $\pp$.

\subsection{Proof of Theorem \ref{thm:MainW}} \label{sect:ProofMainDoob}

To verify (\ref{eq:MainDoobFin}), start from the remarks on the definition of $\ud_{[\pms]}^{(W)}$ at the end of Subsection \ref{sect:metric}, to write
\begin{eqnarray}
\ud_{[\pms]}^{(W)}\big{(} q(\xisun), \delta_{\empiric} \big{)} &:=& \sum_{k=1}^{\infty} \frac{1}{2^k \|h_k\|_{BL}} \Big{|} \int_{\pms} h_k(\pfrak) q(\xisun, \ud\pfrak)  - h_k(\empiric)\Big{|} \leq \int_{\pms} \ud_{\pms}^{(W)}\big{(} \pfrak, \empiric \big{)} q(\xisun, \ud\pfrak) \nonumber \\
&=& \ee\left[\ud_{\pms}^{(W)}\big{(} \randommeasure, \empiric \big{)}\ |\ \xisun\right] \nonumber
\end{eqnarray}
where the inequality holds since the functions $h_k$ are Lipschitz-continuous and the last equality follows from the so-called disintegration theorem (see, e.g., Theorem 6.4 in Kallenberg (2002)). At this stage, one proves that, for every $n_0 > e^e$, $S_{n_0}^{\ast} := \sup_{n \geq n_0} \sqrt{\frac{n}{\log\log n}} \ud_{\pms}^{(W)}\big{(} \randommeasure, \empiric \big{)}$ has finite expectation. In fact, looking at the expectation of $S_{n_0}^{\ast}$ as expectation of the conditional expectation of $S_{n_0}^{\ast}$ given $\randommeasure$, one can take advantage of the fact that the $\tilde{\xi}_i$'s are conditionally i.i.d. given $\randommeasure$ with common p.d. $\randommeasure$. Hence,
\begin{eqnarray}
\ee[S_{n_0}^{\ast}\ |\ \randommeasure] &\leq& \sum_{k=1}^{+\infty} \frac{1}{2^k} \ee\left[\sup_{n \geq n_0} \frac{\big{|} \sum_{i=1}^n [g_k^{\ast}(\tilde{\xi}_i) - \int_{\samplespace} g_k^{\ast}(x) \randommeasure(\ud x)]\big{|}}{\sqrt{n \log\log n}}\ |\ \randommeasure\right] \nonumber \\
&\leq& \sum_{k=1}^{+\infty} \frac{1}{2^k} \left(\ee\left[\sup_{n \geq n_0} \frac{\big{|} \sum_{i=1}^n [g_k^{\ast}(\tilde{\xi}_i) - \int_{\samplespace} g_k^{\ast}(x) \randommeasure(\ud x)]\big{|}^r}{(n \log\log n)^{r/2}}\ |\ \randommeasure\right]\right)^{1/r} \nonumber
\end{eqnarray}
hold for every $r > 2$ and, in view of the theorem recalled in Subsection \ref{sect:teicher}, the last term turns out to be less than a non-random constant, with probability one. Obviously, the same constant represents an upper bound also for $\ee[S_{n_0}^{\ast}]$. This paves the way for mimicking the same argument as in Blackwell and Dubins (1962), to obtain
\begin{eqnarray}
\limsup_{n \rightarrow +\infty} \sqrt{\frac{n}{\log \log n}} \ee\left[\ud_{\pms}^{(W)}\big{(} \randommeasure, \empiric \big{)}\ |\ \xisun\right] &\leq& \lim_{k \rightarrow +\infty} \sup_{n \geq k, l \geq k} \ee\left[\sqrt{\frac{n}{\log \log n}} \ud_{\pms}^{(W)}\big{(} \randommeasure, \empiric \big{)}\ |\ \tilde{\boldsymbol{\xi}}^{(l)}\right] \nonumber \\
&\leq& \lim_{k \rightarrow +\infty} \sup_{l \geq k} \ee\left[\sup_{n \geq k} \sqrt{\frac{n}{\log \log n}} \ud_{\pms}^{(W)}\big{(} \randommeasure, \empiric \big{)}\ |\ \tilde{\boldsymbol{\xi}}^{(l)}\right] \nonumber \\
&\leq& \lim_{k \rightarrow +\infty} \sup_{l \geq k} \ee\left[S_{n_0}^{\ast}\ |\ \tilde{\boldsymbol{\xi}}^{(l)}\right] \nonumber
\end{eqnarray}
the last inequality being valid for every $n_0 > e^e$. Combination of a well-known L\'{e}vy martingale convergence theorem with the $\sigma(\xib)$-measurability of $\randommeasure$ (see (\ref{eq:dF2})) yields $\limsup_{k \rightarrow +\infty} \ee\left[S_{n_0}^{\ast}\ |\ \tilde{\boldsymbol{\xi}}^{(k)}\right] = S_{n_0}^{\ast}$ and then
$$
\limsup_{n \rightarrow +\infty} \sqrt{\frac{n}{\log \log n}}  \ud_{[\pms]}^{(W)}\big{(} q(\xisun), \delta_{\empiric} \big{)}
\leq \limsup_{n_0 \rightarrow +\infty} S_{n_0}^{\ast} \leq \sum_{k=1}^{+\infty} \frac{1}{2^k} \limsup_{n \rightarrow +\infty}
G_{k, n}^{\ast}
$$
with $G_{k, n}^{\ast} := \frac{\big{|} \sum_{i=1}^n [g_k^{\ast}(\tilde{\xi}_i) - \int_{\samplespace} g_k^{\ast}(x) \randommeasure(\ud x)]\big{|}}{\sqrt{n \log\log n}}$. To conclude, notice that
$$
\sigma_k(\randommeasure) := \left(\int_{\samplespace} \left[g_k^{\ast}(x) - \int_{\samplespace} g_k^{\ast}(x) \randommeasure(\ud x)\right]^2 \randommeasure(\ud x)\right)^{1/2} \leq 1
$$
and $\rho\{\lim_{n \rightarrow +\infty} G_{k, n}^{\ast} \leq \sqrt{2} \sigma_k(\randommeasure)\} = \ee[\rho\{\lim_{n \rightarrow +\infty} G_{k, n}^{\ast} \leq \sqrt{2} \sigma_k(\randommeasure)\ |\ \randommeasure\}]$, where, from the Hartman-Wintner law of iterated logarithm, $\rho\{\lim_{n \rightarrow +\infty} G_{k, n}^{\ast} \leq \sqrt{2} \sigma_k(\randommeasure)\ |\ \randommeasure\} = 1$ with probability one, which establishes (\ref{eq:MainDoobFin}).

As to the proof of (\ref{eq:MainDoobPred}), write
\begin{eqnarray}
\ud_{[\samplespace^m]}^{(W)}\big{(} p_m(\xisun), \empiric^m \big{)} &:=& \sum_{k_1, \dots, k_m \in \naturals} \frac{1}{2^{k_1 + \dots + k_m}} \Big{|} \int_{\samplespace^m} \left[\prod_{i=1}^m g_{k_i}^{\ast}(x_i)\right] p_m(\xisun, \ud x_1 \dots \ud x_m)
\nonumber \\
&-& \int_{\samplespace^m} \left[\prod_{i=1}^m g_{k_i}^{\ast}(x_i)\right] \empiric(\ud x_1) \dots \empiric(\ud x_m)\Big{|} \nonumber \\
&=& \sum_{k_1, \dots, k_m \in \naturals} \frac{1}{2^{k_1 + \dots + k_m}} \Big{|} \ee\left[\prod_{i=1}^m g_{k_i}^{\ast}(\tilde{\xi}_{n+i})\ |\ \xisun \right] - \prod_{i=1}^m \left(\frac{1}{n} \sum_{j=1}^n g_{k_i}^{\ast}(\tilde{\xi}_j)\right) \Big{|} \nonumber
\end{eqnarray}
and notice that
$$
\ee\left[\prod_{i=1}^m g_{k_i}^{\ast}(\tilde{\xi}_{n+i})\ |\ \xisun \right] = \ee\left\{ \ee\left[\prod_{i=1}^m g_{k_i}^{\ast}(\tilde{\xi}_{n+i})\ |\ \randommeasure, \xisun \right]\ |\ \xisun \right\} = \ee\left\{ \prod_{i=1}^m
\ee\left[g_{k_i}^{\ast}(\tilde{\xi}_{n+1})\ |\ \randommeasure \right]\ |\ \xisun \right\}
$$
where the last equality follows from the fact that the $\tilde{\xi}_i$'s are conditionally i.i.d., given $\randommeasure$, through Proposition 6.6 in Kallenberg (2002). Now, Lemma 1 of Section 27 in Billingsley (1995) entails
\begin{eqnarray}
\ud_{[\samplespace^m]}^{(W)}\big{(} p_m(\xisun), \empiric^m \big{)} &\leq& \sum_{k_1, \dots, k_m \in \naturals} \frac{1}{2^{k_1 + \dots + k_m}} \ee\left\{ \sum_{i=1}^m \Big{|} \ee\left[g_{k_i}^{\ast}(\tilde{\xi}_{n+1})\ |\ \randommeasure \right] - \frac{1}{n} \sum_{j=1}^n g_{k_i}^{\ast}(\tilde{\xi}_j) \Big{|}\ |\ \xisun \right\} \nonumber \\
&=& \sum_{i=1}^m \ee\left\{ \sum_{k = 1}^{+\infty} \frac{1}{2^k} \Big{|} \ee\left[g_k^{\ast}(\tilde{\xi}_{n+1})\ |\ \randommeasure \right] - \frac{1}{n} \sum_{j=1}^n g_k^{\ast}(\tilde{\xi}_j) \Big{|}\ |\ \xisun \right\} \nonumber \\
&=& m \ee\left[\ud_{\pms}^{(W)}\big{(} \randommeasure, \empiric \big{)}\ |\ \xisun\right] \nonumber
\end{eqnarray}
and the proof of (\ref{eq:MainDoobPred}) can be carried out, from here on, exactly in the same way as the proof of (\ref{eq:MainDoobFin}).

\subsection{Proof of Theorem \ref{thm:MainG}} \label{sect:ProofMainGini}

To prove (\ref{eq:MainGiniFin}), start from the Kantorovich-Rubinstein theorem (see, e.g., Section 11.8 of Dudley (2002))
to write
$$
\ud_{[\pms]}^{(G_1)}\big{(} q(\xisun), \delta_{\empiric} \big{)} = \sup_{\substack{h : \pms \rightarrow \reals \\ \|h\|_{Lip} \leq 1}} \Big{|} \int_{\pms} h(\pfrak) q(\xisun, \ud\pfrak)  - h(\empiric)\Big{|} \leq \ee\left[\ud_{\pms}^{(G_1)}\big{(} \randommeasure, \empiric \big{)}\ |\ \xisun\right]\ .
$$
Now, from a well-known theorem by Dall'Aglio (see, e.g., Chapter 5 in Rachev et. al. (2013)) valid when $\samplespace = \reals$, one has
$$
\ud_{\pms}^{(G_1)}\big{(} \randommeasure, \empiric \big{)} = \int_{\reals} \big{|} \tilde{\textsf{F}}(x) - \frac{1}{n} \sum_{j=1}^n \ind_{(-\infty; x]}(\tilde{\xi}_j) \big{|} \ud x
$$
where $\tilde{\textsf{F}}(x) := \randommeasure((-\infty; x])$ is a random probability distribution function (p.d.f.). Put
$$
S_{n_0}^{\ast} := \sup_{n \geq n_0} \sqrt{\frac{n}{\log\log n}} \int_{\reals} \big{|} \tilde{\textsf{F}}(x) - \frac{1}{n} \sum_{j=1}^n \ind_{(-\infty; x]}(\tilde{\xi}_j) \big{|} \ud x
$$
and exploit the fact that, conditionally on $\tilde{\textsf{F}}$, the $\tilde{\xi}_i$'s are i.i.d. with common p.d.f. $\tilde{\textsf{F}}$ to establish an upper bound for $\ee(S_{n_0}^{\ast})$, as follows
\begin{eqnarray}
\ee(S_{n_0}^{\ast}) &\leq& \ee\left\{ \int_{\reals} \sup_{n \geq n_0} \sqrt{\frac{n}{\log\log n}} \big{|} \tilde{\textsf{F}}(x) - \frac{1}{n} \sum_{j=1}^n \ind_{(-\infty; x]}(\tilde{\xi}_j) \big{|} \ud x \right\} \nonumber \\
&\leq& \ee\left\{ \int_{\reals} \left( \ee\left[\sup_{n \geq n_0} \frac{\big{|} \sum_{j=1}^n [\ind_{(-\infty; x]}(\tilde{\xi}_j) - \tilde{\textsf{F}}(x)] \big{|}^r}{(n \log\log n)^{r/2}}\ \Big{|}\ \tilde{\textsf{F}} \right] \right)^{1/r} \ud x \right\} \nonumber
\end{eqnarray}
with $r > 2$. From a combination of (\ref{eq:teicherpiu2}) with the disintegration theorem one gets
\begin{eqnarray}
\left( \ee\left[\sup_{n \geq n_0} \frac{\big{|} \sum_{j=1}^n [\ind_{(-\infty; x]}(\tilde{\xi}_j) - \tilde{\textsf{F}}(x)] \big{|}^r}{(n \log\log n)^{r/2}}\ \Big{|}\ \tilde{\textsf{F}} \right]\right)^{1/r} &\leq& \left[\alpha_0 \sigma(\tilde{\textsf{F}}(x))^r + \alpha_1 (\sigma(\tilde{\textsf{F}}(x))^r)^{1- \lceil r \rceil} \mu_r(\tilde{\textsf{F}}(x))^{\lceil r \rceil}\right]^{1/r} \nonumber \\
&\leq& \alpha_0^{1/r} \sigma(\tilde{\textsf{F}}(x)) + \alpha_1^{1/r} \sigma(\tilde{\textsf{F}}(x))^{-2} [\mu_r(\tilde{\textsf{F}}(x))]^{3/r} \label{eq:fretta}
\end{eqnarray}
for every $r \in (2, 3]$, with $\sigma(\textsf{F}) := \sqrt{\textsf{F}(x)[1 - \textsf{F}(x)]}$ and $\mu_r(\textsf{F}(x)) :=
[1 - \textsf{F}(x)]^r \textsf{F}(x) + [\textsf{F}(x)]^r [1 - \textsf{F}(x)]$. In view of these remarks, resuming integration with respect to $x$ yields
\begin{eqnarray}
\int_{\reals} \sqrt{\textsf{F}(x)[1 - \textsf{F}(x)]} \ud x &=& \int_{\reals} \sqrt{\textsf{F}(x)[1 - \textsf{F}(x)]} \frac{\sqrt{1 + |x|^{1+\epsilon}}}{\sqrt{1 + |x|^{1+\epsilon}}} \ud x \nonumber \\
&\leq& \left(\int_{\reals} \frac{1}{1 + |x|^{1+\epsilon}} \ud x\right)^{1/2} \left(\int_{\reals} \textsf{F}(x)[1 - \textsf{F}(x)] (1 + |x|^{1+\epsilon}) \ud x\right)^{1/2}  \nonumber \\
&\leq& 2 \left(\int_{0}^{+\infty} [1 - \textsf{F}(x) + \textsf{F}(-x)] (1 + x^{1+\epsilon}) \ud x\right)^{1/2} \nonumber \\
&=& 2 \left(\int_{\reals} \left[ |x| + \frac{|x|^{2+\epsilon}}{(2+\epsilon)} \right] \ud \textsf{F}(x) \right)^{1/2} \label{eq:Cauchy}
\end{eqnarray}
for every $\epsilon \in (0, \delta \wedge 1]$, thanks to a combination of the Cauchy-Schwartz inequality with a well-known representation of moments as in Lemma 1, Section 6.2 of Chow and Teicher (1997). Moreover, for any $r = 2 + \eta \in (2, 3)$,
a combination of that lemma with the H\"{o}lder inequality gives
\begin{eqnarray}
&& \int_{\reals} (\textsf{F}(x)[1 - \textsf{F}(x)])^{-1} \{[1 - \textsf{F}(x)]^r \textsf{F}(x) + [\textsf{F}(x)]^r [1 - \textsf{F}(x)]\}^{3/r} \ud x \leq 2^{\frac{3}{r}} \int_{\reals} [1 - \textsf{F}(x)]^{\frac{3-r}{r}} [\textsf{F}(x)]^{\frac{3-r}{r}} \ud x \nonumber \\
&\leq&  \left(\int_{\reals} \left(\frac{1}{1 + |x|^{1+\epsilon}}\right)^{(1-\eta)/(1+\eta)} \ud x\right)^{(1-\eta)/2} \left(\int_{\reals} \textsf{F}(x)[1 - \textsf{F}(x)] (1 + |x|^{1+\epsilon}) \ud x\right)^{(1+\eta)/2}  \nonumber
\end{eqnarray}
for every $\epsilon \in (0, \delta \wedge 1]$ and $\eta \in (0, 1)$ such that $\frac{1-\eta}{1+\eta}(1+\epsilon) > 1$. Therefore, after bounding the last term in (\ref{eq:fretta}) as above and taking expectation, one arrives at the conclusion that
\begin{equation} \label{eq:Oxford}
\ee\left[\sup_{n > e^e} \sqrt{\frac{n}{\log \log n}}\ud_{\pms}^{(G_1)}\big{(} \randommeasure, \empiric \big{)}\right] < +\infty
\end{equation}
leading to the applicability of the same Blacwell and Dubins argument as in the previous subsection. Whence,
$$
\limsup_{n \rightarrow +\infty} \sqrt{\frac{n}{\log \log n}} \ud_{[\pms]}^{(G_1)}\big{(} q(\xisun), \delta_{\empiric} \big{)} \leq \limsup_{n \rightarrow +\infty} \sqrt{\frac{n}{\log \log n}} \ud_{\pms}^{(G_1)}\big{(} \randommeasure, \empiric \big{)}
$$
and, in view of (\ref{eq:Oxford}), one can combine the extended monotone convergence theorem for decreasing sequences, the law of iterated logarithm and (\ref{eq:Cauchy}) to obtain
\begin{eqnarray}
\limsup_{n \rightarrow +\infty} \sqrt{\frac{n}{\log \log n}} \ud_{\pms}^{(G_1)}\big{(} \randommeasure, \empiric \big{)} &\leq&
\int_{\reals} \limsup_{n \rightarrow +\infty} \sqrt{\frac{n}{\log\log n}} \big{|} \tilde{\textsf{F}}(x) - \frac{1}{n} \sum_{j=1}^n \ind_{(-\infty; x]}(\tilde{\xi}_j) \big{|} \ud x \nonumber \\
&\leq& \sqrt{2} \int_{\reals} \sigma(\tilde{\textsf{F}}(x)) \ud x \leq \left(8 \int_{\reals} \left[ |x| + \frac{|x|^{2+\epsilon}}{(2+\epsilon)} \right] \ud \tilde{\textsf{F}}(x) \right)^{1/2} \nonumber \ .
\end{eqnarray}

To prove (\ref{eq:MainGiniPred}), one can resort to a general argument which shows how to bound $\ud_{[\reals^m_{\sigma}]}^{(G_1)}(p_m(\xisun), \empiric^m)$ in terms of $\ud_{[[\reals]]}^{(G_1)}(q(\xisun), \delta_{\empiric})$. With a view to further developments, to proof will be framed in an abstract setting. At the beginning, fix $m \in \naturals$ and use again the Kantorovich-Rubinstein theorem to write
$$
\ud_{[\samplespace^m_{\sigma}]}^{(G_1)}(p_m(\xisun), \empiric^m) = \sup_{\substack{h : \samplespace^m_{\sigma} \rightarrow \reals \\
h \in Lip_1(\ud_{\samplespace^m_{\sigma}}^{(G_1)})}} \Big{|} \int_{[\samplespace]} \left(\int_{\samplespace^m} h(x) \pfrak^m(\ud x)\right) q(\xisun, \ud\pfrak) - \int_{[\samplespace]} \left(\int_{\samplespace^m} h(x) \pfrak^m(\ud x)\right) \delta_{\empiric}(\ud\pfrak)\Big{|}
$$
where $Lip_1(\ud_{\samplespace^m_{\sigma}}^{(G_1)})$ stands for the class of 1-Lipschitz functions based on the metric $\ud_{\samplespace^m_{\sigma}}^{(G_1)}$ on $\samplespace_{\sigma}^m$. Thus, to carry out the proof, it suffices to verify that, for any fixed $h \in Lip_1(\ud_{\samplespace^m_{\sigma}}^{(G_1)})$, the map $F_h := [\samplespace]_1 \ni \pfrak \mapsto \int_{\samplespace^m_{\sigma}} h(x) \pfrak^m(\ud x)$ is well-defined and Lipschitz-continuous with Lipschitz norm not greater than one. Since the argument is rather technical, its complete explanation will be presented, within a more general framework, in a paper in preparation. Here, one confines oneself to mentioning its basic steps. First, the desired property of $F_h$ is proved, assuming that $\samplespace$ does not contain any isolated point, only for all discrete uniform distributions in $\pms$. In fact, thanks to a Birkhoff theorem on optimal matching (see, e.g., Ambrosio, Gigli and Savar\'{e} (2008) or Villani (2003)), one can write
\begin{eqnarray}
\Big{|} F_h(\frac{1}{N} \sum_{i=1}^N \delta_{x_i}) - F_h(\frac{1}{N} \sum_{i=1}^N \delta_{y_i}) \Big{|}
&\leq& \frac{1}{N^m} \sum_{\substack{i_1, \dots, i_m \\ \in \{1, \dots, N\}}} \ud_{[\samplespace]}^{(G_1)}\left(\frac{1}{m}\sum_{k=1}^m \delta_{x_{i_k}}, \frac{1}{m}\sum_{k=1}^m \delta_{y_{\tau(i_k)}}\right)
\nonumber \\
&\leq& \ud_{[\samplespace]}^{(G_1)}\left(\frac{1}{N} \sum_{i = 1}^N \delta_{x_i}, \frac{1}{N} \sum_{i = 1}^N \delta_{y_i}\right)
\nonumber
\end{eqnarray}
where $\tau$ stands for the optimal coupling permutation. Second, under the same non-isolation extra-condition, the result on
$F_h$ is extended to the whole $[\samplespace]_1$, by means of a suitable density argument. Third, the aforesaid extra-condition is bypassed by reducing the original problem to an equivalent one on $\overline{\samplespace} := \samplespace \times [-1, 1]$
metrized by $\mathrm{D}_{\overline{\samplespace}}((x, s), (y, t)) := \ud_{\samplespace}(x, y) + |s - t|$, to result in a space without isolated points. In fact, extension of any $\mu \in [\samplespace]$ to $\overline{\mu} := \mu \otimes \delta_0 \in [\overline{\samplespace}]$ yields $\ud_{[\overline{\samplespace}]_1}^{(G_1)}(\overline{\mu}_1, \overline{\mu}_2) =
\ud_{[\samplespace]_1}^{(G_1)}(\mu_1, \mu_2)$. Then, one extends $\ud_{\samplespace^m_{\sigma}}^{(G_1)}$ to
$\overline{\ud}_{\samplespace^m_{\sigma}}^{(G_1)}((\xi_1, \dots, \xi_m), (\eta_1, \dots, \eta_m)) := \ud_{[\overline{\samplespace}]}^{(G_1)}\left(\frac{1}{m} \sum_{i=1}^m \delta_{\xi_i}, \frac{1}{m} \sum_{i=1}^m \delta_{\eta_i}\right)$ and any $h : \samplespace_{\sigma}^m \rightarrow \reals$ in $Lip_1(\ud_{\samplespace^m_{\sigma}}^{(G_1)})$ to an $\overline{h} : \overline{\samplespace}_{\sigma}^m \rightarrow \reals$ in
$Lip_1(\overline{\ud}_{\samplespace^m_{\sigma}}^{(G_1)})$ by virtue of Proposition 11.2.3 in Dudley (2002). Finally, the map $F_h$ to $\overline{F}_{\overline{h}}(\overline{\mu}) := \int_{\overline{\samplespace}^m} \overline{h}(x) \overline{\mu}^m(\ud x)$ which satisfies $\big{|} \overline{F}_{\overline{h}}(\overline{\mu}_1) - \overline{F}_{\overline{h}}(\overline{\mu}_2)\big{|} \leq \ud_{[\overline{\samplespace}]}^{(G_1)}(\overline{\mu}_1, \overline{\mu}_2)$ for any $\overline{\mu}_1, \overline{\mu}_2 \in [\overline{\samplespace}]$ and the reasoning is completed by observing that $F_h(\mu) = \overline{F}_{\overline{h}}(\overline{\mu})$ for every $\mu \in [\samplespace]_1$.

\subsection{Proof of Theorem \ref{thm:MainP}} \label{sect:ProofMainProk}

Start by considering the sequence $\{A_{m, j}\}_{j = 1, \dots, k_m + 1}$ of partitions of $\hat{\samplespace}$, where $A_{m, k_m+1} := \mathbb{K}_m^c$, and choose a point $a_{m,j}$ in each set $A_{m,j}$. Since $\mathbb{K}_m^c \downarrow \emptyset$ in view of the $\sigma$-compactness of $\hat{\samplespace}$, one has that, for any $\pfrak \in \pms$, $\pfrak^{(m)} := \sum_{j = 1}^{k_m + 1} \pfrak(A_{m, j} \cap \samplespace) \delta_{a_{m,j}} \Rightarrow \pfrak$ as $m \rightarrow +\infty$. Then, recall the definition of the Fortet-Mourier distance
$$
\ud_{[\Space]}^{(FM)}(\mu_1, \mu_2) := \sup_{\substack{h :\ \Space \rightarrow \reals \\
\| h \|_{BL} \leq 1}} \Big{|} \int_{\Space} h(x) \mu_1(\ud x) - \int_{\Space} h(x) \mu_2(\ud x) \Big{|}\ \ \ \ \ (\mu_1, \mu_2 \in [\Space])\ ,
$$
and exploit the relation $\ud_{[\Space]}^{(P)} \leq [\frac{3}{2} \ud_{[\Space]}^{(FM)}]^{1/2}$ (see Section 11.3 of Dudley (2002)), to obtain
\begin{eqnarray}
\left(\frac{n}{\log n}\right)^{1/4} \ud_{[[\samplespace]]}^{(P)}\left(q(\xisun), \delta_{\empiric} \right) &\leq& \left(\frac{3}{2} \sqrt{\frac{n}{\log n}} \ud_{[[\samplespace]]}^{(FM)}\left(q(\xisun), \delta_{\empiric} \right) \right)^{1/2} \nonumber \\
&\leq& \left(\frac{3}{2} \sqrt{\frac{n}{\log n}} \ee[\ud_{\pms}^{(P)}(\randommeasure, \empiric)\ |\ \xisun] \right)^{1/2} \ . \label{eq:BDProk}
\end{eqnarray}
At this stage, an application of Theorem 11.3.3 in Dudley (2002) shows that
$$
\ud^{(P)}_{\pms}(\randommeasure, \empiric) = \lim_{m \rightarrow +\infty} \ud^{(P)}_{\pms}\left(\sum_{j=1}^{k_m + 1}
\randommeasure(A_{j,m} \cap \samplespace) \delta_{a_{j,m}}, \frac{1}{n}\sum_{i=1}^{n}\sum_{j=1}^{k_m}\ind_{A_{j,m}}(\tilde{\xi}_i) \delta_{a_{j,m}}\right)
$$
and, after putting $\tilde{p}_{j,m} := \randommeasure(A_{j,m} \cap \samplespace)$ and $\tilde{X}_i^{(j,m)} := \ind_{A_{j,m}}(\tilde{\xi}_i)$, an equality displayed on page 95 of Regazzini and Sazonov (2001) yields
$$
\ud^{(P)}_{\pms}\left(\sum_{j=1}^{k_m + 1} \randommeasure(A_{j,m} \cap \samplespace) \delta_{a_{j,m}}, \frac{1}{n}\sum_{i=1}^{n}\sum_{j=1}^{k_m + 1}\ind_{A_{j,m}}(\tilde{\xi}_i) \delta_{a_{j,m}}\right) =
\frac{1}{2} \sum_{j=1}^{k_m + 1} \Big{|}\tilde{p}_{j,m} - \frac{1}{n}\sum_{i=1}^{n} \tilde{X}_i^{(j,m)}\Big{|}\ .
$$
Therefore, since $\sup_{n \geq n_0} \liminf_{m \rightarrow +\infty} x_{n,m} \leq \liminf_{m \rightarrow +\infty} \sup_{n \geq n_0} x_{n,m}$ holds for any subset $\{x_{n,m}\}_{n, m \in \naturals}$ of the real numbers, one gets
\begin{eqnarray}
\sup_{n \geq n_0} \sqrt{\frac{n}{\log n}} \ud_{\pms}^{(P)}(\randommeasure, \empiric) &\leq& \frac{1}{2} \liminf_{m \rightarrow +\infty} \sup_{n \geq n_0} \sqrt{\frac{n}{\log n}} \sum_{j=1}^{k_m + 1} \Big{|}\tilde{p}_{j,m} - \frac{1}{n}\sum_{i=1}^{n} \tilde{X}_i^{(j,m)}\Big{|} \nonumber \\
&\leq& \frac{1}{2} \liminf_{m \rightarrow +\infty} \sum_{j=1}^{k_m + 1} \left( \sup_{n \geq n_0} \sqrt{\frac{n}{\log n}}\ \Big{|}\tilde{p}_{j,m} - \frac{1}{n}\sum_{i=1}^{n} \tilde{X}_i^{(j,m)}\Big{|} \right)\ . \nonumber
\end{eqnarray}
Now, combine the Lyapunov inequality for moments with (\ref{eq:newteicherpiu2}) to write, for any fixed $r > 2$,
\begin{eqnarray}
&& \ee\left[ \sup_{n \geq n_0} \sqrt{\frac{n}{\log n}}\ \Big{|}\tilde{p}_{j,m} - \frac{1}{n}\sum_{i=1}^{n} \tilde{X}_i^{(j,m)}\Big{|}\ |\ \randommeasure \right] \nonumber \\
&\leq& \left(\ee\left[ \sup_{n \geq n_0} \frac{\Big{|}\sum_{i=1}^{n}(\tilde{X}_i^{(j,m)} - \tilde{p}_{j,m}) \Big{|}^r}{(n \log n)^{r/2}} \ \big{|}\ \randommeasure \right]\right)^{1/r} \nonumber \\
&\leq& \gamma(r, n_0) \left[\sqrt{\tilde{p}_{j,m}(1 - \tilde{p}_{j,m})} + [\tilde{p}_{j,m}^r(1 - \tilde{p}_{j,m}) + \tilde{p}_{j,m}(1 - \tilde{p}_{j,m})^r]^{1/r}\right] \nonumber
\end{eqnarray}
with a suitable non-random constant $\gamma(r, n_0)$. Since
$$
\liminf_{m \rightarrow +\infty} \sum_{j=1}^{k_m + 1} \left[\sqrt{\tilde{p}_{j,m}(1 - \tilde{p}_{j,m})} + [\tilde{p}_{j,m}^r(1 - \tilde{p}_{j,m}) + \tilde{p}_{j,m}(1 - \tilde{p}_{j,m})^r]^{1/r}\right] \leq 3Y_r(\randommeasure)
$$
holds true, an application of the conditional Fatou lemma shows that
$$
\ee\left[\sup_{n \geq n_0} \sqrt{\frac{n}{\log n}} \ud_{\pms}^{(P)}(\randommeasure, \empiric)\right] \leq \frac{3}{2} \gamma(r, n_0) \ee[\Pi_r(\randommeasure)] < +\infty\ .
$$
Consequently, one can resort to the already utilized Blackwell-Dubins argument to obtain (\ref{eq:MainProkFin}) directly from (\ref{eq:BDProk}).

As to (\ref{eq:MainProkPred}), one gets its validity directly from (\ref{eq:MainProkFin}). Indeed, since it has already proved at the end of Subsection \ref{sect:ProofMainGini} that $\ud_{[\samplespace^m_{\sigma}]}^{(G_1)}(p_m(\xisun), \empiric^m) \leq
\ud_{[\pms]}^{(G_1)}(q(\xisun), \delta_{\empiric})$ is in force for every $m \in \naturals$, one can simply resort to the bounds
$\ud_{[\Space]}^{(P)} \leq [\frac{3}{2} \ud_{[\Space]}^{(FM)}]^{1/2} \leq [\frac{3}{2} \ud_{[\Space]}^{(G_1)}]^{1/2}$ and
$$
\ud_{[[\samplespace]]}^{(G_1)}\left(q(\xisun), \delta_{\empiric} \right) \leq \ee[\ud_{\pms}^{(P)}(\randommeasure, \empiric)\ |\ \xisun] \ .
$$

\subsection{Proof of Proposition \ref{prop:PropP}} \label{sect:sobolev}

If $\randommeasure(A) = \tilde{\lambda} \sum_{i = 1}^{\tilde{N}} \tilde{p}_i \delta_{\tilde{x}_i}(A) + (1 - \tilde{\lambda}) \int_A \tilde{f}(x) \ud x$ for every $A \in \mathscr{B}(\samplespace)$, one can simply exploit the concavity of the function $x \mapsto x^{1/r}$, $x \in [0, +\infty)$, to get
$$
\Pi_r(\randommeasure) \leq [\tilde{\lambda}]^{1/r} \liminf_{m \rightarrow +\infty} \sum_{j = 1}^{k_m + 1} \sum_{i = 1}^{\tilde{N}} [\tilde{p}_i]^{1/r} \delta_{\tilde{x}_i}(A_{m, j})\ + [1 - \tilde{\lambda}]^{1/r} \liminf_{m \rightarrow +\infty} \sum_{j = 1}^{k_m + 1} \left(\int_{A_{m, j}} [(\tilde{f}(x))^{1/r}]^r \ud x\right)^{1/r} \ .
$$
Since $\sum_{j = 1}^{k_m + 1} \delta_{\tilde{x}_i}(A_{m, j}) = 1$ for all $i \in \{1, \dots, \tilde{N}\}$ and $m \in \naturals$, it is enough to study the second summand in the right-hand side of the above inequality. Therefore, if $l := \frac{r-1}{r}d$ is an integer, it is easy to show that $\frac{1}{r} = 1 - \frac{l}{d}$, which, in conjunction with the regularity assumptions on $\samplespace$, guarantees the validity of the Sobolev imbedding $\mathrm{W}^{l, 1}(\samplespace) \subset \mathrm{L}^r(\samplespace)$. See Adams and Fournier (2003) for more information. To conclude, upon noting that the imbedding constants can be fixed independently of the partitions, it is enough to observe that
$$
\sum_{j = 1}^{k_m + 1} \int_{A_{m,j}} |D^{\boldsymbol{\alpha}} (\tilde{f}(x))^{1/r}| \ud x =
\int_{\samplespace} |D^{\boldsymbol{\alpha}} (\tilde{f}(x))^{1/r}| \ud x
$$
holds for every multi-index $\boldsymbol{\alpha}$.

\vspace{1cm}

\noindent \textbf{Acknowledgements}
We thank Pietro Rigo and Giuseppe Savar\'e for helpful discussions. Work partially supported by MIUR-2008MK3AFZ and INdAM-GNAMPA Project 2015.

\vspace{2cm}

           \textsc{Donato Michele Cifarelli} \\
           Dipartimento di Scienze delle Decisioni \\
           Universit\`a commerciale L. Bocconi, 20136 Milano, Italy \\
           \\
           \textsc{Emanuele Dolera} \\
           Dipartimento di Scienze Fisiche, Informatiche e Matematiche \\
           Universit\`a di Modena e Reggio Emilia, 41125 Modena, Italy \\
           \emph{emanuele.dolera@unimore.it, emanuele.dolera@unipv.it}             \\
           \\
           \textsc{Eugenio Regazzini} \\
           Dipartimento di Matematica \\
           Universit\`a di Pavia, 27100 Pavia, Italy   \\
           Affiliated also with CNR-IMATI, Milano, Italy \\
           \emph{eugenio.regazzini@unipv.it}


\begin{thebibliography}{99}

\bibitem{} \textsc{Adams, R.A.} and \textsc{Fournier, J.J.F.} (2003). \emph{Sobolev spaces}. $2^{nd}$ ed. Academic Press, Amsterdam.

\bibitem{} \textsc{Aldous, D.J.} (1985). Exchangeability and related topics. \emph{\'{E}cole d'\'{e}t\'{e} de probabilit\'{e}s de Saint-Flour, XIII-1983}. Lecture Notes in Math. $\mathbf{1117}$ 1-198. Springer, Berlin.

\bibitem{} \textsc{Ambrosio, L.}, \textsc{Gigli, N.} and \textsc{Savar\'{e}, G.} (2008). \emph{Gradient flows in Metric Spaces and in the Space of Probability Measures}. $2^{nd}$ ed. Birkh\"{a}user, Basel.

\bibitem{} \textsc{Barron, A., Schervish, M.J.} and \textsc{Wasserman, L.} (1999). The consistency of posterior distribution in nonparametric problems. \emph{Ann. Statist}. $\mathbf{27}$ 536-561

\bibitem{} \textsc{Bassetti, F.} (2011). Quantitative comparisons between finitary posterior distributions and Bayesian posterior distributions. \emph{J. Statist. Plann. Inference} $\mathbf{141}$ 787-799

\bibitem{} \textsc{Bernstein, S.N.} (1917). \emph{Theory of Probability} (Russian). Moscow.

\bibitem{} \textsc{Berti, P., Crimaldi, I., Pratelli, L.} and \textsc{Rigo, P.} (2009). Rate of convergence of predictive distributions for dependent data. \emph{Bernoulli} $\mathbf{15}$ 1351-1367

\bibitem{} \textsc{Billingsley, P.} (1995). \emph{Probability and Measure}. $3^{rd}$ ed. Wiley, New York.

\bibitem{} \textsc{Billingsley, P.} (1999). \emph{Convergence of Probability Measures}. $2^{nd}$ ed. Wiley, New York.

\bibitem{} \textsc{Blackwell, D.} and \textsc{Dubins, L.E.} (1962). Merging of opinions with increasing information.
\emph{Ann. Math. Statist.} $\mathbf{33}$ 882-886

\bibitem{} \textsc{Chow, Y.S.} and \textsc{Teicher, H.} (1997). \emph{Probability Theory. Independence, Interchangeability, Martingales}. $3^{rd}$ ed. Springer, New York.

\bibitem{} \textsc{Daley, D.J.} and \textsc{Vere-Jones, D.} (2003). \emph{An Introduction to the Theory of Point Processes} Vol. 1. $2^{nd}$ ed. Springer, New York.

\bibitem{} \textsc{Diaconis, P.} and \textsc{Freedman, D.} (1986). On consistency of Bayes estimates. \emph{Ann. Statist.} $\mathbf{14}$, 1-26

\bibitem{} \textsc{Diaconis, P.} and \textsc{Freedman, D.A.} (1990). On the uniform consistency of Bayes estimates for multinomial probabilities \emph{Ann. Statist}. $\mathbf{18}$ 1317-1327

\bibitem{} \textsc{Doob, J.L.} (1949). Application of the theory of martingales. In \emph{Le calcul des Probabilit\'es et ses Applications} 23-27. Colloques Internationaux du Centre National de la Recherche Scientifique, Paris.

\bibitem{} \textsc{Dudley, R.M.} (1969). The speed of mean Glivenko-Cantelli convergence. \emph{Ann. Math. Statist.} $\mathbf{40}$ 40-50

\bibitem{} \textsc{Dudley, R.M.} (2002). \emph{Real Analysis and Probability}. Cambridge University Press, Cambridge.

\bibitem{} \textsc{Efron, B.} (2003). Robbins, empirical Bayes and microarrays. \emph{Ann. Statist}. $\mathbf{31}$ 366-378

\bibitem{} de \textsc{Finetti, B.} (1929). Funzione caratteristica di un fenomeno aleatorio (Italian). \emph{Atti del Congresso Internazionale dei Matematici}, Bologna September $3^{rd}-10^{th}$ 1928, 179-190

\bibitem{} de \textsc{Finetti, B.} (1930). Funzione caratteristica di un fenomeno aleatorio (Italian). \emph{Atti Reale Accademia Nazionale dei Lincei, Mem.} $\textbf{4}$ 86-133

\bibitem{} de \textsc{Finetti, B.} (1933). Sull'approssimazione empirica di una legge di probabilit\`a (Italian). \emph{Giornale dell'Istituto Italiano degli Attuari} $\textbf{4}$ 415-420

\bibitem{} de \textsc{Finetti, B.} (1937). La pr\'{e}vision: ses lois logiques, ses sources subjectives (French). \emph{Ann. Inst. H. Poincar\'{e}} $\textbf{7}$ 1-68

\bibitem{} de \textsc{Finetti, B.} (1952). La notion de ``distribution d'opinions'' comme base d'un essai d'interpretation de la Statistique (French). \emph{Publ. Inst. Stat. Univ. Paris} \textbf{1} 1-19

\bibitem{} de \textsc{Finetti, B.} (1970). \emph{Teoria delle Probabilit\`a}. Einaudi, Torino. Two volumes [English translation, \emph{Theory of Probability} (1975). Wiley, New York. Two volumes.]

\bibitem{} \textsc{Freedman, D.} (1999). Wald Lecture: On the Bernstein-von Mises theorem with infinite-dimensional parameter. \emph{Ann. Statist}. $\mathbf{27}$ 1119-1141

\bibitem{} \textsc{Ghosal, S., Ghosh, J.K.} and \textsc{van der Vaart, A.W.} (2000). Convergence rates of posterior distributions. \emph{Ann. Statist}. $\mathbf{28}$ 500-531

\bibitem{} \textsc{Ghosh, J.K.} and \textsc{Ramamoorthi, R.V.} (2003). \emph{Bayesian Nonparametrics}. Springer, New York.

\bibitem{} \textsc{Gini, C.} (1914) Di una misura della dissomiglianza tra due gruppi di quantit\`a a delle sue applicazioni allo studio delle relazioni statistiche (Italian). Atti del R. Istit. Veneto di Scienze, Lettere e Arti. Tomo LXXIV,

\bibitem{} \textsc{Kallenberg, O.} (2002). \emph{Foundations of Modern Probability}. $2^{nd}$ ed. Springer-Verlag, New York.

\bibitem{} \textsc{Laplace, P.S.} (1812). \emph{Th\'{e}orie Analytique des Probabilit\'{e}s} (French). Courcier, Paris.

\bibitem{} \textsc{Leti, G.} (1961). Nuovi tipi di distanze fra insiemi di punti e loro applicazioni alla statistica (Italian). \emph{Metron} $\textbf{21}$ 139-169

\bibitem{} \textsc{Leti, G.} (1962). Il termine generico delle tabelle di cograduazione e di contrograduazione (Italian). Biblioteca del \emph{Metron} $\textbf{C 1}$ 253-277.

\bibitem{} \textsc{Lo, A.Y.} (1984). On a class of Bayesian nonparametric estimate: I. Density estimates. \emph{Ann. Statist.} $\mathbf{12}$, 351-357

\bibitem{} von \textsc{Mises, R.} (1919). Grundlagen der Wahrscheinlichkeitsrechnung. \emph{Math. Z}.  $\mathbf{5}$ 52-100

\bibitem{} von \textsc{Mises, R.} (1964). \emph{Probability and Statistics. General}. American Mathematical Society, Providence.

\bibitem{} \textsc{Petrov, V.V.} (1995). \emph{Limit Theorems of Probability Theory. Sequences of Independent Random Variables}. Oxford University Press, New York.

\bibitem{} \textsc{Poincar\'e, H.} (1912). \emph{Calcul des Probabilit\'{e}s}. Gauthier-Villars, Paris.

\bibitem{} \textsc{Rachev, S.T., Klebanov, L.B., Stoyanov, S.V.} and \textsc{Fabbozzi, F.J.} (2013). \emph{The Methods of Distances in the Theory of Probability and Statistics}. Springer, New York.

\bibitem{} \textsc{Regazzini, E.} and \textsc{Sazonov, V.V.} (2001). Approximation of laws of random probabilities by mixtures of Dirichlet distributions with applications to nonparametric bayesian inference. \emph{Theory Probab. Appl.} $\mathbf{45}$ 93-110

\bibitem{} \textsc{Robbins, H.} (1956). An empirical Bayes approach to statistics. \emph{Proc. Third Berkeley Symp. Math. Statist. Probab}. $\mathbf{1}$ 157-163

\bibitem{} \textsc{Robbins, H.} (1964). The empirical Bayes approach to statistical decision problems. \emph{Ann. Math. Statist}. $\mathbf{35}$ 1-20.

\bibitem{} \textsc{Romanovsky, V.} (1931). Sulle probabilit\`a ``a posteriori'' (Italian). \emph{Giornale dell'Istituto Italiano degli Attuari}. $\mathbf{2}$ 493-511

\bibitem{} \textsc{Schwartz, L.} (1965). On Bayes procedures. \emph{Z. Wahrsch. Verw. Gebiete}. $\mathbf{4}$ 10-26

\bibitem{} \textsc{Stroock, D.W.} (2001). \emph{Probability theory: an analytical view}. $2^{nd}$ ed. Cambridge Universuty Press, Cambridge.

\bibitem{} van der \textsc{Vaart, A.} and \textsc{Wellner, J.A.} (1996). \emph{Weak convergence and Empirical Processes}. Springer, New York.

\bibitem{} \textsc{Villani, C.} (2003). \emph{Topics in Optimal Transportation}. American Mathematical Society, Providence.

\bibitem{} \textsc{Yukich, J.E.} (1989). Optimal matching and empirical measures. \emph{Proc. Amer. Math. Soc.} $\mathbf{107}$ 1051-1059.
\end{thebibliography}
\end{document}